\def\qed{{~~~\vrule height .75em width .4em}}
\def\irr#1{{\rm Irr}(#1)}
\def\cent#1#2{{\bf C}_{#1}(#2)}
\def\syl#1#2{{\rm Syl}_#1(#2)}
\def\nor{\triangleleft\,}
\def\fit#1{{\bf F}(#1)}
\def\oh#1#2{{\bf O}_{#1}(#2)}
\def\Oh#1#2{{\bf O}^{#1}(#2)}
\def\zent#1{{\bf Z}(#1)}
\def\ker#1{{\rm ker}(#1)}
\def\cd#1{{\rm cd}(#1)}
\def\dl#1{{\rm dl}(#1)}
\def\norm#1#2{{\bf N}_{#1}(#2)}
\def\ritem#1{\item{{\rm #1}}}
\def\iitem#1{\goodbreak\par\noindent{\bf #1}}
\def\aut#1{{\rm Aut}(#1)}
\def\sbs{\subseteq}
\def\dimm#1#2{{\rm dim}_{#1}(#2)}
\def\sps{\supseteq}
\let\phi=\varphi

\def\irrs#1#2{{\rm Irr}({#2}|{#1})}
\def\cds#1#2{{\rm cd}({#2}|{#1})}
\def\van#1{{\bf V}(#1)}
\def\AA{{\cal A}}
\def\BB{{\cal B}}
\def\ref#1{{\bf [#1]}}
\def\gar{3}
\def\book{6}
\def\berk{2}
\def\berg{1}
\def\gluck{4}
\def\isa{5}

\magnification = \magstep1

{\nopagenumbers
\vglue 1.5truein
\font\bb=cmbx12

\centerline{{\bb IRREDUCIBLE CHARACTER DEGREES}}
\smallskip
\centerline{{\bb AND NORMAL SUBGROUPS}}
\medskip
\centerline{by}
\bigskip
\centerline{{\bf I.\ M. Isaacs and Greg Knutson}}
\bigskip\bigskip

\centerline{Mathematics Department}
\centerline{University of Wisconsin}
\centerline{480 Lincoln Drive}
\centerline{Madison WI~~~53706}
\centerline{USA}
\smallskip
\centerline{E-mail:~~isaacs@math.wisc.edu}
\vfil
\font\sf=cmr9
\sf
The research of the first author was supported by a grant from the
United States NSA.
\eject}

\iitem{1. Introduction.}

Let $G$ be a finite group and, as usual, write $\cd G$ to denote the
set of degrees of the irreducible characters of $G$. This set of positive
integers encodes a great deal of information about the structure of $G$,
and we mention just a few of the many known results of this type. 
If $G$ is nilpotent (or more generally, if $G$ is an M-group), a result of
K.~Taketa asserts that the derived length $\dl G$ is at most equal to
$|\cd G|$. (See Theorem~5.12 of \ref\book.) It has been conjectured
that the Taketa inequality $\dl G \le |\cd G|$ holds for all solvable
groups, but this remains unproved. (It is known, however, that
that $\dl G \le 2|\cd G|$ for all solvable groups $G$ \ref\gluck,
and the Taketa inequality is known to hold if $|G|$ is odd \ref\berg.)

It is an old result of the first author that if $|\cd G| \le 3$, then $G$ is
necessarily solvable and that $\dl G \le |\cd G|$. (See Theorems~12.15
and~12.6 of \ref\book.) The Taketa inequality is also known to hold
for solvable groups $G$ if $|\cd G| = 4$. (This result of S.~Garrison is
the principal theorem in his Ph.D.\ thesis \ref\gar, which remains
unpublished.)

If more information about the set $\cd G$ is known than its cardinality,
then, of course, one can expect to be able to deduce correspondingly
more  information about $G$. For example, J.~Thompson proved that
if there is some prime $p$ that divides every member of $\cd G$
exceeding $1$, then $G$ has a normal $p$-complement.
(See Corollary~12.2 of \ref\book.)

Our goal in this paper is to study how the structure of a normal subgroup
of $G$ is influenced by the degrees of an appropriate subset of $\irr G$.
It seems reasonable that the characters that should be relevant to
controlling the structure of $N \nor G$ are exactly those whose kernels
do not contain $N$, and so we introduce some convenient notation.
Given that $N \nor G$, we write
$$
\irrs NG = \{\chi \in \irr G | N \not\sbs \ker\chi\}
$$
and
$$
\cds NG = \{\chi(1) | \chi \in \irrs NG\} \,.
$$

If we take $N = G'$, the derived subgroup of $G$, we see that
$\irrs NG$ is exactly the set of nonlinear characters in $\irr G$, and
thus $|\cds NG| = |\cd N| - 1$ in this case. Also, if $G$ is solvable and
$N = G'$, then $\dl N = \dl G - 1$, and hence the conjectured Taketa
inequality, that $\dl G \le |\cd G|$ for solvable groups, can be restated as
$\dl N \le |\cds NG|$, where, of course, $N = G'$.

The first of our main results provides an upper bound for $\dl N$ in terms
of $|\cds NG|$, where $N$ is an arbitrary solvable normal subgroup of
an arbitrary finite group $G$. We do not need to assume that $N = G'$
or that $G$ is solvable.
\medskip

\iitem{THEOREM A.}~~{\sl Let $N \nor G$, where $N$ is solvable,
and write $n = |\cds NG|$. Then $|\dl N| \le f(n)$ for some quadratic
function $f$. Furthermore, if $G$ is assumed to be solvable, then a
linear upper bound exists.}
\medskip

Returning to the case where $N = G'$, we know that if $|\cds NG| \le 2$,
which corresponds to the situation where $|\cd G| \le 3$, then $N$ is
solvable and $\dl N \le \cds NG$. We get exactly this result for
arbitrary normal subgroups $N$ of arbitrary finite groups $G$.
\medskip

\iitem{THEOREM B.}~~{\sl Let $N \nor G$ and suppose that
$|\cds NG| \le 2$. Then $N$ is solvable and $\dl N \le |\cds NG|$.}
\medskip

Unfortunately, our proof of the solvability of $N$ in the case where
$|\cds NG| = 2$ in Theorem~B relies on a consequence the
classification of simple groups.

If we assume some solvability hypothesis on the whole group $G$,
we can push this one step farther.
\medskip

\iitem{THEOREM C.}~~{\sl Let $N \nor G$ and assume that $G$ is
$p$-solvable for all prime divisors $p$ of $|N|$. If $|\cds NG| = 3$, then
$\dl N \le 3$.}
\medskip

Note that in the situation of Theorem~C, the subgroup $N$ is
automatically solvable. If the whole group $G$ is solvable and
$|\cd G| = 4$, we can apply Theorem~C to the subgroup $N = G'$,
and we deduce that $\dl G \le 4$. In other words, Theorem~C includes
the unpublished result of Garrison's thesis as a special case. It seems
reasonable to conjecture that Theorem~C might remain true if the
hypothesis that $G$ is $p$-solvable for all prime divisors $p$ of
$|N|$ were replaced by the weaker assumption that $N$ is solvable,
but we have been unable to prove this. Perhaps it is even true
when $N$ is solvable that $\dl N \le |\cds NG|$, even when $|\cds NG|$,
exceeds $3$.

There is at least one result already in the literature that provides
control over the structure of a subgroup $N \nor G$ in terms of
$\cds NG$. This result of Y.~Berkovich \ref\berk, is the analog of the
theorem of Thompson that we mentioned previously. Our results are
dependent on Berkovich's theorem, and so for completeness, we
present a simple proof of a somewhat strengthened version.
\medskip

\iitem{THEOREM D. (Berkovich)}~~{\sl Let $N \nor G$ and suppose
that every member of $\cds{N'}G$ is divisible by some fixed prime $p$.
Then $N$ has a normal $p$-complement.}
\medskip

Since $\irrs{G'}G$ is exactly the set of nonlinear irreducible characters
of $G$, we see that in the case $N = G$, Theorem~D is exactly
Thompson's theorem. What Berkovich actually proved was that
the intersection of the kernels of the nonlinear irreducible
characters of $G$ having $p'$-degree has a normal $p$-complement.
In other words, Berkovich showed that $N \nor G$ has a
normal $p$-complement if all nonlinear members of $\irrs NG$ have
degree divisible by $p$. Theorem~D is stronger, however,
because the hypothesis is weaker: our assumption is that irreducible
characters of $G$ that lie over nonlinear irreducible characters of $N$
should have degree divisible by $p$.

We close this introduction with a brief explanation of how this paper
came to be written. At the time of his death in 1996, Greg Knutson was
Isaacs' Ph.D.\ student at the University of Wisconsin, Madison. Knutson
had made substantial progress in his study of the set $\cds NG$, and it
seemed appropriate for Isaacs to continue that research and to write
what might have become Knutson's thesis. The present paper is the
result of that effort. 
\bigskip

\iitem{2. Berkovich's Theorem.}

In this section, we give a quick proof of our slightly generalized
version of Berkovich's theorem. \medskip

\iitem{Proof of Theorem D.}~~Write $M = \Oh pN$, let
$P \in \syl p M$ and choose $S \in \syl pG$ such that $S \sps P$. Our
goal is to show that $p$ does not divide $|M|$, and so we assume that
$P > 1$ and we work to obtain a contradiction. Note that
$P = S \cap M  \nor S$, and thus $[P,S] < P$, and we can choose a
nonprincipal $S$-invariant linear character $\lambda$ of $P$.

Now $S$ stabilizes $\lambda^M$, and thus it permutes the irreducible
constituents of this character. But $\lambda^M(1) = |M:P|$ is not
divisible by $p$, and so $\lambda^M$ must have some $S$-invariant
constituent $\alpha \in \irr M$ with $\alpha(1)$ not divisible by $p$.
Now $\alpha$ is invariant in $MS$, and the $p$-power $|MS:M|$ is
is relatively prime to $\alpha(1)$. Furthermore, the determinantal order
$o(\alpha)$ is not divisible by $p$ since $M = \Oh pM$, and thus
$o(\alpha)$ is also relatively prime to $|MS:M|$. We deduce
that $\alpha$ extends to some character $\beta \in \irr{MS}$.
(See Corollary~6.28 of \ref\book.)

Next, observe that $\beta^G(1) = \beta(1)|G:MS| = \alpha(1)|G:MS|$
is not divisible by $p$, and hence it has some constituent $\chi \in \irr G$
with degree not divisible by $p$. By hypothesis, $\chi$ is not a member
of $\irrs{N'}G$, and thus $N' \sbs \ker\chi$, and so the irreducible
constituents of $\chi_N$ are linear. The irreducible constituents of
$\chi_M$ are therefore linear, and in particular, $\alpha$ is linear. Thus
$\alpha$ is an extension of $\lambda$ to $M$, and it follows that
$o(\lambda)$ divides $o(\alpha)$, which, as we have seen, is not
divisible by $p$, and therefore $o(\lambda)$ is not divisible by $p$.
This is a contradiction since $\lambda$ is a nontrivial linear character
of a $p$-group, and this completes the proof.\qed
\bigskip

\iitem{3. An analog of Taketa's theorem.}

The principal result of this section is the following theorem, whose proof
is a variation on the standard Taketa argument used to show that
$\dl G \le |\cd G|$ when $G$ is an M-group. To state the result, we shall
use the notation $m_p$ to denote the $p$-part of a positive integer $m$,
where $p$ is a prime number. In other words, $m_p$ is the largest
power of $p$ that divides $m$.
\medskip

\iitem{(3.1) THEOREM.}~~{\sl Let $N \nor G$, where $N$ has an
abelian normal $p$-complement for some prime $p$. If $\chi \in \irrs NG$
is chosen such that $\chi(1)_p$ is as small as possible, then
$N' \sbs \ker\chi$. In particular, if $N > 1$, then $\cds{N'}G < \cds NG$.}
\medskip

\iitem{Proof.}~~Let $P/N \in \syl p{G/N}$ and choose an irreducible
constituent $\psi$ of $\chi_P$ with $\psi(1)$ as small as possible.
By hypothesis, $N$ has an abelian normal $p$-complement, and since
$P/N$ is a $p$-group, we see that $P$ also has an abelian normal
$p$-complement, and it follows that $P$ is an M-group. (See
Theorems~6.22 and~6.23 of \ref\book, for example.) We can thus write
$\psi = \lambda^P$, where $\lambda$ is a linear character of some
subgroup $Q \sbs P$.

If $\eta$ is an irreducible constituent of the character $(1_Q)^P$, we
claim that $N \sbs \ker\eta$. We can certainly assume that $\eta$ is
nonprincipal, and thus $\eta(1) < |P:Q| = \psi(1)$. Also, since $P$ has
an abelian normal $p$-complement, Ito's theorem guarantees that
all irreducible characters of $P$ have $p$-power degrees. (See
Theorem~6.15 of \ref\book.) In particular, $\eta(1)$ is a $p$-power, and
since $p$ does not divide $|G:P|$, we conclude that $\eta(1)$ is exactly
the $p$-part of $\eta^G(1) = \eta(1)|G:P|$. It follows that there exists an
irreducible constituent $\xi$ of $\eta^G$ such that $\xi(1)_p \le \eta(1)$.

As all irreducible constituents of $\chi_P$ have $p$-power degrees and
$\psi(1)$ is the smallest of these degrees, we see that $\psi(1)$ divides
$\chi(1)$. Therefore, $\xi(1)_p \le \eta(1) < \psi(1) \le \chi(1)_p$, and it
follows from the choice of $\chi$ that $\xi \not\in \irrs NG$, and hence
$N \sbs \ker\xi$. But $\eta$ is an irreducible constituent of $\xi_P$, and
thus $N \sbs \ker\eta$, as claimed. Since this holds for every
irreducible constituent $\eta$ of $(1_Q)^P$, we deduce that
$N \sbs \ker{(1_Q)^P} \sbs Q$, and thus $N' \sbs Q' \sbs \ker\lambda$.

Since $N' \nor G$, it follows that $N'$ is contained in the kernel of every
irreducible constituent of $\lambda^G = \psi^G$. In particular,
$N' \sbs \ker\chi$, as required.

Finally, given that $N > 1$, we see that $\cds NG$ is nonempty. But
$\cds{N'}G \sbs \cds NG$, and this subset is missing all members of
$\cds NG$ that have the smallest possible $p$-part. This completes
the proof.\qed \medskip

We can now prove the case of Theorem~B where $|\cds NG| \le 1$.
\medskip

\iitem{(3.2) COROLLARY.}~~{\sl Let $N \nor G$ and suppose that
$|\cds NG| \le 1$. Then $\dl N \le |\cds NG|$, and in particular,
$N$ is abelian.}
\medskip

\iitem{Proof.}~~If $N = 1$, then $\dl N = 0$ and the desired inequality
holds. We thus assume that $N > 1$, and since $|\cds NG| = 1$ in this
case, it suffices to show that $N$ is abelian, and we do this by induction
on $|N|$.  If $N$ is nonabelian, it has a nonlinear irreducible
character $\psi$, and we choose a prime divisor $p$ of $\psi(1)$.
If $\chi \in \irr G$ lies over $\psi$, then $N \not\sbs \ker\chi$, and thus
$\chi \in \irrs NG$ and $\chi(1)$ is the unique member of $\cds NG$.
Also, $\psi(1)$ divides $\chi(1)$, and so $p$ divides every member of
$\cds NG$. By Theorem~D, therefore, $N$ has a normal
$p$-complement $M$. Also, $\psi(1)$ divides $|N|$, and so $p$
divides $|N|$, and $M < N$.

Clearly, $\irrs MG \sbs \irrs NG$, and thus $|\cds MG| \le 1$ and the
inductive hypothesis guarantees that $M$ is abelian. Thus $N$ has
an abelian normal $p$-complement, and hence $\cds{N'}G < \cds NG$
by Theorem~3.1. It follows that $\cds{N'}G$ is empty, and therefore
$N' = 1$ and $N$ is abelian. This is a contradiction, and the theorem
is proved.\qed.
\medskip

We also have another easy consequence of Theorem~3.1.
\medskip

\iitem{(3.3) COROLLARY.}~~{\sl Let $N \nor G$ and suppose that
$N$ is nilpotent. Then $\dl N  \le |\cds NG|$.}
\medskip

\iitem{Proof.}~~Since $N$ is nilpotent, its derived length is the
maximum of the derived lengths of its Sylow subgroups. Also, if
$P \in \syl pN$, then $P \nor G$ and $\cds PG \sbs \cds NG$. It suffices,
therefore, to prove the theorem in the case where $N$ is a $p$-group,
and we do this by induction on $|N|$. 

Assuming that $N$ is a $p$-group and that $N > 1$, we see from
Theorem~3.1 that $\cds{N'}G < \cds NG$, and thus
$|\cds{N'}G| < |\cds NG|$. Also, $N' < N$, and so by the inductive
hypothesis, we have $\dl{N'} \le |\cds{N'}G| \le |\cds NG| - 1$, and it
follows that $\dl N \le |\cds NG|$, as required.\qed
\bigskip

\iitem{4. The Fitting subgroup.}

It is known that for any nonidentity solvable group $G$, the set
$\cd{G/F}$ fails to contain the largest member of $\cd G$, where
$F = \fit G$ is the Fitting subgroup. In this section, we prove the
analogous result for $\cds NG$.
\medskip

\iitem{(4.1) THEOREM.}~~{\sl Let $N \nor G$ and choose
$\chi \in \irrs NG$ such that $\chi(1) = \max(\cds NG)$. 
Then $N \cap \ker\chi$ is nilpotent. Also, if $K/M$ is an abelian chief
factor of $G$, where $M \sbs N \cap \ker\chi$, then $K$ is nilpotent.}
\medskip

The proof we present for Theorem~4.1 is a minor variation on the
argument used to prove Theorem~12.19 of \ref\book, which depends
on the `vanishing-off' subgroup $\van\chi$ of a character $\chi \in \irr G$.
Recall that by definition, $\van\chi$ is the subgroup generated by all
elements $x \in G$ such that $\chi(x) \ne 0$, and note that it is obvious
that $\ker\chi \sbs \van\chi$. We shall need the following variation on
Lemma~12.18 of \ref\book.
\medskip

\iitem{(4.2) LEMMA.}~~{\sl Let $V = \van\chi$, where $\chi \in \irr G$.
Suppose that $K \nor G$ and that $K \cap V \sbs \ker\chi$. Then
$K \sbs \ker\chi$.}
\medskip

\iitem{Proof.}~~Since $\chi$ vanishes on $K - V$, we see that
$|K|[\chi_K,1_K] = |K \cap V|[\chi_{K \cap V}, 1_{K \cap V}] > 0$, where
the inequality holds since $K \cap V \sbs \ker\chi$. It follows that
$\chi_K$ has a principal constituent, and thus
$K \sbs \ker\chi$.\qed
\medskip

\iitem{Proof of Theorem 4.1.}~~Assuming that the theorem is false, we
define normal subgroups $K$ and $M$ of $G$ as follows. If
$N \cap \ker\chi$ is not nilpotent, define $K = N \cap \ker\chi = M$.
Otherwise, let $K/M$ be an abelian chief factor of $G$, where
$M \sbs N \cap \ker\chi$ and $K$ is not nilpotent. Note that our
definitions imply that although $K$ is not nilpotent, $M$ is nilpotent
in the case where $M < K$.

Since $K$ is not nilpotent, we can choose a nonnormal Sylow
$p$-subgroup $P$ of $K$ and we write $H = \norm GP$.
If $K = M$, then by the Frattini argument, we have
$G = MH$. If, on the other hand, $K > M$, then $M$ is nilpotent and
$K/M$ is an abelian chief factor of $G$. In that situation, $P \not\sbs M$,
and thus $p$ must be the unique prime divisor of $K/M$, and
$K = MP$. By the Frattini argument, we deduce that
$G = KH = MPH = MH$, and thus $G = MH$ in all cases.

Since $G = MH$ and $M \sbs \ker\chi$, it follows that $\chi_H$ is
irreducible, and we write $\psi = \chi_H$. By definition, $M \sbs N$,
and thus $N = M(N \cap H)$. Since $N \not\sbs \ker\chi$ and
$M \sbs \ker\chi$, we see that $N \cap H \not\sbs \ker\chi$, and thus
$N \cap H \not\sbs \ker\psi$. It follows that $N \cap H$ is not contained
in the kernel of any irreducible constituent of $\psi^G$, and in particular,
$N$ is not contained in any of these kernels. All irreducible constituents
of $\psi^G$, therefore, lie in $\irrs NG$, and hence each of them
has degree at most $\max(\cds NG) = \chi(1) = \psi(1)$. It follows that
the restriction to $H$ of each irreducible constituent of $\psi^G$ must be
$\psi$ exactly, and thus by Lemma~12.17 of \ref\book, we conclude that
$\van\psi$ is normal in $G$. We shall show that
$K \cap H \sbs \van\psi$, and thus $K \cap H = K \cap \van\psi \nor G$.
This will yield a contradiction, however, since $K \cap H = \norm KP$ is
not normal in $K$.

Consider the subgroup $X = M(K \cap \van\psi)$, which is normal
in $G$. Note that $M \sbs X \sbs K$, so that either $X = K$ or
$X = M$, and we suppose first that $X = K$. Then since
$K \cap \van\psi \sbs H$, we have
$$
K \cap H = M(K \cap \van\psi) \cap H =
(M \cap H)(K \cap \van\psi) \sbs \van\psi  \,,
$$
where the containment holds because
$M \cap H \sbs \ker\psi \sbs \van\psi$. We are done in this case,
and so we can assume that $X = M$, so that $K \cap \van\psi \sbs M$.

We now have $(K \cap H) \cap \van\psi \sbs M \cap H \sbs \ker\psi$, 
and we can apply Lemma~4.2 to the normal subgroup $K \cap H$
in the group $H$. We deduce that
$K \cap H \sbs \ker\psi \sbs \van\psi$, as wanted. The proof is now
complete.\qed
\medskip

\iitem{(4.3) COROLLARY.}~~{\sl Let $F = \fit N$, where $1 < N \nor G$
and $N$ is solvable, and let $m = \max(\cds NG)$. Then
$m \not\in \cds{(N/F)}{(G/F)}$, which is therefore a proper subset of
$\cds NG$.}
\medskip

\iitem{Proof.}~~Note that $m$ exists since $N > 1$.
If we view irreducible characters of $G/F$ as irreducible characters
of $G$, we see that $\irrs{(N/F)}{(G/F)} \sbs \irrs NG$, and thus
$\cds{(N/F)}{(G/F)} \sbs \cds NG$. It therefore suffices to show
that if $\chi \in \irrs NG$ with $\chi(1) = m$, then $F \not\sbs \ker\chi$.

Suppose that $F \sbs \ker\chi$. Since $\chi \in \irrs NG$, we know that
$N \not\sbs \ker\chi$, and thus $F < N$, and we can choose a chief
factor $K/F$ of $G$ with $K \sbs N$. Since $N$ is solvable, $K/F$ is
abelian, and thus Theorem~4.1 applies and $K$ is nilpotent. This is a
contradiction, however, since $K \nor N$, but $K \not\sbs F = \fit N$.
This completes the proof.\qed
\medskip

Our next result includes part of Theorem~A. If $N$ is solvable, we bound
both the Fitting height $h(N)$ and the derived length $\dl N$ in terms of
$|\cds NG|$.
\medskip

\iitem{(4.4) COROLLARY.}~~{\sl Let $N \nor G$, where
$N$ is solvable, and write $n = |\cds NG|$. Then $h(N) \le n$ and
$\dl N \le n(n+1)/2$.}
\medskip

\iitem{Proof.}~~We proceed by induction on $|N|$. If $N = 1$, then
$h(N) = 0 = \dl N$ and there is nothing to prove. Otherwise, set 
$F = \fit N$ and observe that $F > 1$. Then $|\cds{(N/F)}{(G/F)}| \le n-1$
by Corollary~4.3, and we have
$$
h(N) = 1 + h(N/F) \le 1 + (n-1) = n \, ,
$$
by the inductive hypothesis. Also, $\irrs FG \sbs \irrs NG$, and so
$|\cds FG| \le n$ and $\dl F \le n$ by Corollary~3.3. We know from the
inductive hypothesis that $\dl{N/F} \le (n-1)n/2$, and it follows that
$\dl N \le n + (n-1)n/2 = n(n+1)/2$, as required.\qed
\medskip

To complete the proof of Theorem~A, we use a known argument,
independent of the rest of this paper.
\medskip

\iitem{(4.5) THEOREM.}~~{\sl Let $G$ be solvable and suppose
$N \nor G$ and $|\cds NG| = n$. Then $\dl N \le 3n$.}
\medskip

\iitem{Proof.}~~Working by induction on $|N|$, it suffices to show that
$|\cds MG| \le n-1$, where $M = N'''$. Certainly, $\cds MG \sbs \cds NG$,
and so it suffices to show that if $\chi \in \irrs NG$ with
$\chi(1) = \min(\cds NG)$, then $M \sbs \ker\chi$, and thus
$\chi \not\in \irrs MG$. But $N \sbs \ker\psi$ for all $\psi \in \irr G$ with
$\psi(1) < \chi(1)$, and thus Theorem~6 of \ref\isa\ applies. This result
says exactly that $M = N''' \sbs \ker\chi$, as required.\qed
\medskip

We mention that the solvability of $G$ is used in the proof of
Theorem~6 of \ref\isa\ so that the Fong-Swan theorem can be applied
to lift a modular irreducible representation of $G$ to the complex
numbers.

We close this section with an easy, but technical, application of
Theorem~4.1 that will be needed in what follows.
\medskip

\iitem{(4.6) LEMMA.}~~{\sl Let $m = \max(\cds NG)$, where $N \nor G$.
Suppose $M \sbs N$, where $M \nor G$ and $\fit M$ is a $p$-group.
If $m \in \cds{(N/M')}{(G/M')}$, then $M$ is a $p$-group.}
\medskip

\iitem{Proof.}~~We can assume that $N > 1$, so that $m$ exists. By
hypothesis, there exists $\chi \in \irrs NG$ with $\chi(1) = m$ and
$M' \sbs \ker\chi$. Then $M'$ is nilpotent by Theorem~4.1,  and hence
$M' \sbs \fit M$ is a $p$-group. If some prime $q \ne p$ divides $|M/M'|$,
then there exists a subgroup $K \sbs M$ such that $K/M'$ is an abelian
$q$-group that is a chief factor of $G$. We deduce from Theorem~4.1
that $K$ is nilpotent, and thus $K \sbs \fit M$ is a $p$-group, and this is a
contradiction since $|K/M'|$ is divisible by $q \ne p$. It follows that
$|M/M'|$ is not divisible by any prime different from $p$, and hence $M$
is a $p$-group, as required.\qed
\bigskip
 
\iitem{5. Two degrees.}

In this section we study the situation where $|\cds NG| = 2$. We begin
with two easy general lemmas.
\medskip

\iitem{(5.1) LEMMA.}~~{\sl Let $N \nor G$ and suppose that
$a = \min(\cds NG)$. If $H \sbs G$ and $|G:H| \le a$, then
$N \sbs H$. If, in fact, $|G:H| < a$, then $N \sbs H'$.}
\medskip

\iitem{Proof.}~~We can certainly assume that $H < G$, so that
$(1_H)^G$ is reducible. If $\eta$ is any irreducible constituent of
$(1_H)^G$, therefore, we have $\eta(1) < |G:H| \le a$, and thus
$\eta \not\in \irrs NG$. It follows that $N \sbs \ker\eta$, and since $\eta$
was arbitrary, we have $N \sbs \ker{(1_H)^G} \sbs H$, as required.
Finally, if $N \not\sbs H'$, then $N$ is not contained in the kernel of
some linear character $\lambda$ of $H$, and thus the irreducible
constituents of $\lambda^H$ lie in $\irrs NG$. In this case, we have
$a \le \lambda^G(1) = |G:H|$, and the last assertion of the lemma now
follows.\qed
\medskip

\iitem{(5.2) LEMMA.}~~{\sl Let $\alpha \in \irr N$, where $N \nor G$.
Suppose that $\alpha$ is invariant in $G$ and let $p$ be a prime
that does not divide $o(\alpha)\alpha(1)$. Then there exists
$\gamma \in \irr G$, lying over $\alpha$, such that $\gamma(1)$ is not
divisible by $p$.}
\medskip

\iitem{Proof.}~~Let $P/N \in \syl p{G/N}$. By Corollary~6.28 of \ref\book,
we can choose an extension $\beta$ of $\alpha$ to $P$, and we observe
that $\beta^G(1) = \beta(1)|G:P|$ is not divisible by $p$. The degree of
some irreducible constituent $\gamma$ of $\beta^G$ is therefore not
divisible by $p$, and the proof is complete since $\gamma$ lies over
$\alpha$.\qed
\medskip

\iitem{(5.3) THEOREM.}~~{\sl Let $N \nor G$ and suppose that
$|\cds NG| = 2$. If $N$ is solvable, then $\dl N \le 2$.}
\medskip

\iitem{Proof.}~~Assume that $\dl N \ge 3$ and work by induction on
$|N|$. By Corollary~3.3, we see that $N$ cannot be nilpotent, and thus
$F < N$, where we have written $F = \fit N$. If $M \sbs N$ with
$M \nor G$, then $\cds{(N/M)}{(G/M)} \sbs \cds NG$, and hence if
$M >1$, the inductive hypothesis implies that $N/M$ is metabelian,
and so $N'' \sbs M$. In other words, $N''$ is the unique minimal normal
subgroup of $G$ contained in $N$, and it follows that $F$ is a
$p$-group, where $p$ is the unique prime divisor of $|N''|$. 

Now let $M < N$, where $M \nor G$. Since $\cds MG \sbs \cds NG$,
the inductive hypothesis guarantees that $M'$ is abelian, and thus
$N/M'$ must be nonabelian. It follows by Corollary~3.2 that
$|\cds{(N/M')}{(G/M')}| \ge 2$, and thus both members of $\cds NG$ lie
in $\cds{(N/M')}{(G/M')}$. In particular, if we write $\cds NG = \{a,b\}$
with $a < b$, we see that $b \in \cds{(N/M')}{(G/M')}$. Since
$\fit M \sbs \fit N$ is a $p$-group, it follows by Lemma~4.6 that $M$ is a
$p$-group, and thus $M \sbs F$. In other words, $F$ is the unique
normal subgroup of $G$ that is maximal with the property of being
properly contained in $N$. In particular, $N' \sbs F$, and hence $F$ is
nonabelian. Also, $N/F$ is a chief factor of $G$, and hence it is a
$q$-group for some prime $q \ne p$. It follows that $\Oh pN = N$.

Since $F$ is nonabelian, we have $|\cds FG| \ge 2$, and hence there
exists $\psi \in \irrs FG$ with $\psi(1) = a$. Let $\alpha$ be an irreducible
constituent of $\psi_F$, so that $\alpha$ is nonprincipal, and let $T$
be the stabilizer of $\alpha$ in $G$. Then $|G:T|$ divides $\psi(1) = a$,
and hence $T \sps N$ by Lemma~5.1, and thus $\alpha$ is invariant in
$N$. Also, since $F$ is a $p$-group and $N/F$ is a $q$-group with
$q \ne p$, it follows by Corollary~6.28 of \ref\book\ that $\alpha$
extends to $N$ and that a unique extension $\beta$ can be chosen so
that $o(\beta) = o(\alpha)$. But $o(\beta)$ is not divisible by $p$ since
$N = \Oh pN$, and thus $\alpha$ is a nonprincipal irreducible
character of the $p$-group $F$ such that $o(\alpha)$ is not divisible by
$p$. Therefore $\alpha$ cannot be linear, and so $p$ divides
$\alpha(1)$, and thus $p$ divides $a$. Also, since $\beta$ is uniquely
determined by $\alpha$, it is fixed by every element of $G$ that fixes
$\alpha$.

Now let $\lambda$ be a nontrivial linear character of $N$ with
$F \sbs \ker\lambda$, and let $S$ be the stabilizer of $\lambda$ in $G$.
Note that $p$ divides neither $\lambda(1) = 1$ nor $o(\lambda)$, which
is a $q$-power, and so by Lemma~5.2, there exists $\mu \in \irr S$ lying
over $\lambda$ and such that $p$ does not divide $\mu(1)$. Let
$\chi = \mu^G$ and note that $\chi \in \irr G$ by the Clifford
correspondence. 

As $\chi$ lies over $\lambda \in \irr N$ and $\lambda$ is nonprincipal, we
see that $N \not\sbs \ker\chi$ and $\chi \in \irrs NG$. It is not possible,
however, that $\chi(1) = b$ since otherwise, Theorem~4.1 would
imply that $N$ is nilpotent, which is not the case. (Theorem~4.1 applies
because $F \sbs \ker\lambda \sbs \ker\chi$ and $N/F$ is an abelian chief
factor of $G$.) It follows that $|G:S|\mu(1) = \chi(1) = a$, and thus
$|G:S|$ is divisible by the full $p$-part of $a$.

Now consider the character $\lambda\beta \in \irr N$, and observe that
the characters $\beta$ and $\lambda\beta$ uniquely determine
$\lambda$ by Gallagher's theorem. (See Corollary~6.17 of \ref\book).
Let $U$ be the stabilizer of $\lambda\beta$ in $G$ and note that $U$
stabilizes $\alpha = (\lambda\beta)_F$, and thus $U$ also stabilizes
$\beta$. It follows that $U$ stabilizes $\lambda$, and thus $U \sbs S$,
and we conclude that $|G:U|$ is divisible by the full $p$-part of $a$.

Let $\delta \in \irr U$ lie over $\lambda\beta$ and note that $\delta^G$
is irreducible, and in fact, $\delta^G \in \irrs NG$ since $\lambda\beta$ is
nonprincipal. Also, $\delta^G(1) = \delta(1)|G:U|$, and the $p$-part
of this integer exceeds the $p$-part of $a$ because $\delta(1)$ is a
multiple of $\alpha(1)$, which is divisible by $p$. It follows that
$\delta^G(1) \ne a$, and we conclude that $\delta^G(1) = b$, and hence
$p$ divides $b$. Thus $a$ and $b$ are each divisible by $p$, and
therefore $N$ has a normal $p$-complement by Theorem~D. This is
a contradiction, however, since $N = \Oh pN$ and yet $p$ divides $|N|$,
and this completes the proof.\qed
\medskip

To complete the proof of Theorem~B, we need to be able to prove
Theorem~5.3 without assuming that $N$ is solvable. In order to do that,
we need the following result, which is an easy consequence of the
classification of simple groups. 
\medskip

\iitem{(5.4) LEMMA.}~~{\sl Let $G$ be a nonabelian simple group and
suppose that $A$ acts on $G$ via automorphisms and that
$(|A|,|G|) = 1$. Then $A$ has nontrivial fixed points in $G$.}
\medskip

\iitem{Proof.}~~By the classification, the only possibilities for $A$ and
$G$ are that (up to conjugacy in $\aut G$) the action of $A$ is via
field automorphisms on a group of Lie type defined over some field
$E$. In this case, the fixed point subgroup of $A$ in $G$ is exactly the
corresponding group of Lie type defined over the fixed field $F$ of $A$
acting on $E$. In particular, this fixed-point subgroup is nontrivial.\qed
\medskip

\iitem{Proof of Theorem~B.}~~We have $N \nor G$ and $|\cds NG| \le 2$,
and our task is to show that $N$ is solvable and that $\dl N \le |\cds NG|$.
By Theorem~3.2, we can assume that $|\cds NG| = 2$, and we write
$\cds NG = \{a,b\}$ with $a < b$. By Theorem~5.3, we can assume that
$N$ is not solvable.

We proceed by induction on $|N|$. Suppose $1 < M < N$ with
$M \nor G$. Then each of $\cds MG$ and $\cds{(N/M)}{(G/M)}$ is
contained in $\cds NG$, and hence by the inductive hypothesis, each
of $M$ and $N/M$ is solvable. Since we are assuming that $N$ is not
solvable, this situation does not arise, and thus $N$ is a minimal normal
subgroup of $G$, and hence it is a direct product of isomorphic
nonabelian simple groups. In particular, $N$ does not have a normal
$p$-complement for any prime divisor $p$ of $|N|$, and hence by
Theorem~D, no prime divisor of $|N|$ can divide both $a$ and $b$.

Let $\AA \sbs \irr N$ be the set of irreducible constituents of characters
$\chi_N$, where $\chi \in \irrs NG$ has degree $a$. Similarly, let $\BB$
be the set of irreducible characters of $N$ that lie under characters in
$\irrs NG$ that have degree $b$. Then $\AA$ and $\BB$ are nonempty
and $\AA \cup \BB = \irr N - \{1_N\}$. Also, if $\alpha \in \AA$ and
$\beta \in \BB$, then $\alpha(1)$ and $\beta(1)$ must be relatively
prime since any common prime divisor would be a divisor of
$|N|$ that divides both $a$ and $b$, and we know that no such prime
exists. Since $N = N'$, the members of $\AA$ and $\BB$ are all
nonlinear, and it follows that $\AA$ and $\BB$ are disjoint. Also, if $N$ is
a direct product of two or more simple groups, we can find
$\alpha \in \AA$ and $\beta \in \BB$ such that $\alpha\beta$ is
irreducible. Since the degree of this product character is not coprime
either to $\alpha(1)$ or to $\beta(1)$, we see that $\alpha\beta$ cannot
lie in $\AA \cup \BB$, and this is a contribution. We conclude that $N$ is
actually simple.

Suppose that $r$ is a common prime divisor of $a$ and $b$, and choose
$R \in \syl r G$. Now $r$ does not divide $|N|$, and so $(|R|,|N|) = 1$,
and it follows by Lemma~5.4 that $R$ has nontrivial fixed points in $N$.
By the Glauberman correspondence, $R$ fixes some nonprincipal
character $\theta \in \irr N$. (See Theorem~13.1 of \ref\book.) It follows
that $\theta$ extends to $\eta \in \irr{NR}$, and every irreducible
constituent of $\eta^G$ lies in $\irrs NG$. Each of these constituents has
degree $a$ or $b$, and so the degree of each irreducible constituent of
$\eta^G$ is divisible by $r$. It follows that $\eta^G(1) = \theta(1)|G:NR|$
is divisible by $r$, and this is impossible since neither $\theta(1)$ nor
$|G:NR|$ is a multiple of $r$. This contradiction proves that $r$ does not
exist, and thus $a$ and $b$ are relatively prime.

Let $\chi_1$ and $\chi_2$ be characters of degree $a$ in $\irrs NG$
and assume that no constituent of $(\chi_1)_N$ is the complex conjugate
of an irreducible constituent of $(\chi_2)_N$. Then
$(\chi_1\chi_2)_N$ has no principal constituent, and it follows that
$\chi_1\chi_2$ is a sum of characters in $\irrs NG$, with (say) $u$
irreducible constituents of degree $a$ and $v$ irreducible constituents of
degree $b$, counting multiplicities. We thus have
$a^2 = ua + vb$, and we deduce that $a$ divides $v$. But $ab > a^2$,
and thus $v < a$, and we conclude that $v = 0$. This shows that all
irreducible constituents of $\chi_1\chi_2$ have degree $a$. 

Now let $\chi ,\psi \in \irrs NG$ with $\chi(1) = a$ and $\psi(1) = b$.
Since the irreducible constituents of $\chi_N$ and of $\psi_N$ have
different degrees, we see that $(\chi\psi)_N$ can have no principal
constituent, and thus $\chi\psi$ can be written as a sum of (say) $r$
members of $\irrs NG$ of degree $a$ and $s$ of degree $b$,
counting multiplicities. Thus $ab = ra + sb$, and hence $a$ divides $s$
and $b$ divides $r$. It follows from this that either $r = 0$ or $s = 0$.

Suppose $s = 0$. Then $\chi\psi$ is a sum of $b$ irreducible characters
of degree $a$, and we let $\xi$ be any one of these irreducible
constituents. Then $\overline\chi\xi$ has $\psi$ as an irreducible
constituent of degree $b$, and it follows from the earlier calculation that
some irreducible constituent of $(\overline\chi)_N$ is the complex
conjugate of a constituent of $\xi_N$. It follows that $\xi_N$ is of the form
$e\chi_N$ for some rational number $e$, and since this is true for every
choice of $\xi$ as an irreducible constituent of $\chi\psi$, we deduce that
$\chi_N\psi_N = f\chi_N$, where necessarily $f = \psi(1)$. Now $\psi_N$
is a faithful character of $N$, and so its value is different from $f$ on
every nonidentity element of $N$. It follows that $\chi(n) = 0$
for $1 \ne n \in N$, and this is impossible, since $\chi_N$ has no
principal constituent. It follows that $s > 0$ and $r = 0$, and thus
all constituents of $\chi\psi$ have degree $b$. In particular,
if $\alpha \in \AA$ and $\beta \in \BB$, then all irreducible constituents
of $\alpha\beta$ lie in $\BB$.

Now let $\alpha_1,\alpha_2 \in \AA$. If $\alpha_1\alpha_2$ has an
irreducible constituent $\beta \in \BB$, then $\overline{\alpha_1}\beta$
has $\alpha_2$ as a constituent, and that contradicts the result of the
previous paragraph. It follows that the set $\AA \cup \{1_N\} \sbs \irr N$
has the property that the product of any two of its members is a sum of
members of this set. Since $\AA$ is nonempty, we can choose
$\alpha \in \AA$ and it follows that all irreducible constituents of
$\alpha^n$ lie in $\AA \cup \{1_N\}$ for all positive integers $n$. Since
$\alpha$ is faithful, the Burnside-Brauer theorem guarantees that all
irreducible characters of $N$ lie in this set. (See Theorem~4.3 of
\ref\book.) This is a contradiction since $\BB$ is nonempty and disjoint
from $\AA$.\qed
\bigskip

\iitem{6. The exceptional case.}

Let $\chi \in \irrs NG$, where $N \nor G$. Suppose that $X/Y$ is an
abelian chief factor of $G$, where $X \sbs N$, and let $p$ be the prime
divisor of $|X/Y|$. In this situation, we say that $X/Y$ is a
{\bf reducing $p$-section} for $\chi$ in $N$ if the irreducible constituents
of $\chi_Y$ have smaller degrees than the irreducible constituents of
$\chi_X$. (Of course, all irreducible constituents of each of these
restrictions have equal degrees.) As we shall see, if $G$
is $p$-solvable and $\chi(1) = \min(\cds NG)$, then as a general rule, a
reducing $p$-section for $\chi$ in $N$ will be central in $N$. This is not
always true, however, and our goal in this section is to deduce a number
of properties of the exceptional case, where $X/Y$ is not central. It will
turn out, for example, that no exception can occur unless $p = 2$ and
$\chi(1)$ is even.

We begin with an example. Let $G = GL(2,3)$ and take $N = SL(2,3)$,
so that $\cds NG = \{2,3,4\}$. Let $\chi \in \irr G$ be faithful of degree
$2$ and take $X$ to be the quaternion normal subgroup of $N$ of order
$8$ and let $Y = \zent X$, so that $X/Y$ is a chief factor of order $4$ of
$G$. Now $\chi_X$ is irreducible and $\chi_Y$ reduces, and thus $X/Y$
is a reducing $2$-section for $\chi$ in $N$. But $N$ acts nontrivially on
$X/Y$, and so we are in the exceptional case.
\medskip

\iitem{(6.1) THEOREM.}~~{\sl Let $N \nor G$, where $G$ is
$p$-solvable, and let $\chi \in \irrs NG$ with $\chi(1) = \min(\cds NG)$.
Assume that $X/Y$ is a reducing $p$-section for $\chi$ in $N$ but that
$X/Y$ is not central in $N$. The following then occur.
\smallskip
\ritem{(a)} $p = 2$.
\smallskip
\ritem{(b)} $\chi(1)$ is even.
\smallskip
\ritem{(c)} $|\cds NG| \ge 3$.
\smallskip
\ritem{(d)} If $|\cds NG| = 3$, then $\cds NG$ consists of two powers of
$2$ and an odd number. In this case, $N/X$ is abelian.
\medskip}

\iitem{Proof.}~~Since $X/Y$ is an elementary abelian $p$-group,
we can write $|X/Y| = p^e$, where $e$ is a positive integer, and we can
view $X/Y$ as an irreducible $G$-module of dimension $e$ over the
field $F = GF(p)$. In fact, if $E$ is the field of $G$-endomorphisms of
$X/Y$, then $E$ is an extension field of $F$ and $X/Y$ can be viewed
as a vector space over $E$ with dimension $e/|E:F| \le e$. This
makes $X/Y$ into an absolutely irreducible module for $G$ over $E$,
and as such, it corresponds to an irreducible Brauer character $\phi$
of $G$ at the prime $p$, where $\phi(1) = \dimm E{X/Y} \le e$. 

Because $G$ is $p$-solvable, we can apply the Fong-Swan theorem to
find $\alpha \in \irr G$ such that $\alpha(x) = \phi(x)$ for all $p$-regular
elements $x \in G$. In particular, $\alpha(1) = \phi(1) \le e$. We are
assuming that $X/Y$ is a $p$-group that is a chief factor of $G$,
acted on nontrivially by $N \nor G$. It follows that there exists a
$p$-regular element $x \in N$ that acts nontrivially on $X/Y$. Thus
$\alpha(x) = \phi(x) \ne \phi(1) = \alpha(1)$, and we see that
$x \not\in \ker\alpha$. Thus $N \not\sbs \ker\alpha$, and
$\alpha \in \irrs NG$. In particular, if we set $a = \min(\cds NG)$, we see
that $a \le \alpha(1) \le e$.

Now let $\theta$ be an irreducible constituent of $\chi_X$, write
$f = \theta(1)$ and let $t$ be the number of distinct conjugates of
$\theta$ in $G$. Then $a = \chi(1)$ is a multiple of $ft$, and in particular,
$ft \le a$. Also, since $X/Y$ is a reducing section for $\chi$, we see that
$\theta_Y$ reduces, and it follows from the fact that $X/Y$ is a $p$-group
that $\theta$ is induced from a character of some subgroup $U$, with
$Y \sbs U < X$. (This follows, for instance, from Theorem~6.22 of
\ref\book.) Since $U \nor X$, we see that $\theta$ vanishes on $X - U$,
and we let $V$ be the unique smallest subgroup of $X$ containing $Y$
such that $\theta$ vanishes on $X - V$. (Thus $V = Y\van\theta$.) We
now know that $Y \sbs V < X$.

Since $\theta$ uniquely determines $V$, we see that $V$ is normalized
by the stabilizer of $\theta$ in $G$, and thus the number of distinct
conjugates of $V$ is at most $t$. Since $X/Y$ is a chief factor of $G$
and $Y \sbs V < X$, it follows that the intersection of the $G$-conjugates
of $V$ must be exactly $Y$, and therefore, $p^e = |X:Y| \le |X:V|^t$.
We can obtain control over the index $|X:V|$ by observing that since
$\theta$ vanishes on elements of $X - V$, we have
$$
|X| =|X|[\theta,\theta] =
|V|[\theta_V,\theta_V] \le |V|\theta(1)^2 = |V|f^2 \, ,
$$
and so $|X:V| \le f^2$.

Combining the various inequalities we have found, we obtain
$$
p^{ft} \le p^a \le p^e \le |X:V|^t \le f^{2t} \,,
$$
and thus $p^f \le f^2$. Since $2 \le p$, this yields that $2^f \le f^2$, and
thus $f \in \{2,3,4\}$, and if $f$ is equal to either $2$ or $4$, then all of
the previous inequalities are actually equalities and in particular, $p = 2$.
Furthermore, if $f = 3$, then $p^3 \le 3^2$, and in this case too, $p = 2$.
But the latter situation is impossible since $\theta_Y$ is reducible, and
thus $p$ must divide $\theta(1) = f$. We have definitely established,
therefore, that $f \in \{2,4\}$, that $p = 2$ and that all of the equalities
in the previous paragraphs are equalities. In particular, conclusions (a)
and (b) both hold. (Recall that $f$ divides $a = \chi(1)$, and thus
$\chi(1)$ is even.)

Let $T$ be the stabilizer of $\theta$ in $G$, so that $\chi$ is induced
from some character $\xi \in \irr T$, lying over $\theta$. Then
$ft = a = \chi(1) =|G:T|\xi(1) = t\xi(1)$, and so $\xi(1) = f$ and $\xi$ is an
extension of $\theta$ to $T$. Also, $|G:T| = t = a/f < a$, and hence
$N \sbs T'$ by Lemma~5.1. But $X < N$ since $X/Y$ is abelian, but is
not a central factor of $N$, and thus $X < T'$, and we deduce that
$T/X$ is nonabelian. It follows that there exists a nonlinear character
$\beta \in \irr{T/X}$, and we see that $\xi\beta \in \irr T$ lies over $\theta$.
It follows that $(\xi\beta)^G \in \irrs NG$, and thus $\cds NG$ contains
some proper multiple $ra$ of $a$, where $r = \beta(1)$. We see now
that $a$ and $ra$ are two distinct even members of $\cds NG$, and so
to show that $|\cds NG| \ge 3$, it suffices to show that some member of
this set is odd. If not, then $N$ has a normal $2$-complement by
Theorem~D, and this subgroup necessarily centralizes the
$2$-group $X/Y$. This would contradict the previously observed fact
that some $p$-regular element of $N$ must act nontrivially on $X/Y$.
(Recall that $p = 2$.) Thus $\cds NG$ does contain an odd number $m$,
and so $|\cds NG| \ge 3$, proving conclusion (c).

We assume now that $|\cds NG| = 3$, so that we have
$\cds NG = \{a,ra,m\}$. In particular, $r$ is unique, and thus
$\cd{T/X} = \{1,r\}$. In this case $T/X$ must be metabelian, and since
$N \sbs T'$, we deduce that $N/X$ is abelian, as asserted in (d).
To complete the proof, we must show that $a$ and $r$ are each
powers of $2$.

Now $X/V$ is a module for $T/X$ over $GF(2)$, and since $|X/V| = f^2$,
we see that the dimension of this module is $2\log_2(f)$, and this is
equal to $f$ since $f$ is $2$ or $4$. As $N \sbs T$, we see that $\Oh2N$
acts on $X/V$, and we claim that this action is nontrivial. Otherwise,
$[X,\Oh2N] \sbs V$, and thus $[X,\Oh2N] \sbs Y$ because
$[X,\Oh2N] \nor G$ and the intersection of the $G$-conjugates of $V$ is
$Y$. It follows that the $2$-regular elements of $N$ act trivially on $X/Y$,
and we know that this is not the case.

The module $X/V$ for $T/X$ yields a Brauer character $\gamma$
(at the prime $2$), and we know that $\gamma(1) = f$ and that
$\gamma(x) \ne \gamma(1)$ for some $2$-regular element $x \in N$.
We deduce by the Fong-Swan theorem that there exists
$\beta \in \irr{T/X}$ with $N/X \not\sbs \ker\beta$ and $\beta(1) \le f$.
Viewing $\beta \in \irr T$, we see that $N \not\sbs \ker\beta$, and so
each irreducible constituent of $\beta^G$ lies in $\irrs NG$.
We conclude that $a \le \beta^G(1) = |G:T|\beta(1) \le tf = a$. It follows
that $\beta(1) = f$, and we observe also that
$\beta(1) \in \cd{T/X} = \{1,r\}$. Thus $r = f$ is a power of $2$, as desired.

Next, we claim that the product of the $t$ conjugates of $\theta$ is an
irreducible character of $X$. Writing $\delta$ to denote this product, we
observe that $\delta(x) = 0$ unless $x$ is in $V$ and in every
$G$-conjugate of $V$, and thus $\delta$ vanishes on $X - Y$. Also, we
know that $[\theta_V,\theta_V] = \theta(1)^2$, and it follows that
$\theta_V$ is a multiple of a linear character. Therefore, $\theta_Y$ is a
multiple of a linear character, and since a similar conclusion holds for
each conjugate of $\theta$, it follows that $\delta_Y$ is a multiple of a
linear character. Thus $|\delta(y)| = \delta(1)$ for
all elements $y \in Y$, and therefore
$|X|[\delta,\delta] =  |Y|[\delta_Y,\delta_Y] = |Y|\delta(1)^2$. But
$|X:Y| = |X:V|^t = f^{2t} = \delta(1)^2$, and thus
$[\delta,\delta] = 1$ and $\delta$ is irreducible, as claimed.

It follows that $f^t = \delta(1)$ must divide some member of $\cds NG$,
and we conclude that $f^t$ divides $ra = fa = f^2t$. Thus
$f^{t-2} \le t$, and hence $2^{t-2} \le t$ and $2^t \le 4t$. It follows
that $t \le 4$. If $t = 3$, however, then since $f \in \{2,4\}$, we see that
$f^t = f^3$ does not divide $f^2t$, and this is a contradiction. We
conclude that $t \ne 3$, and thus $t \in \{1,2,4\}$ and $a = ft$ is a power
of $2$. This concludes the proof of (d).\qed
\medskip

Note that we actually proved more than we stated in part (d) of
Theorem~6.1. We showed, for example, that the smallest member of
$\cds NG$ is a divisor of $16$ and that the other power of $2$ in this set
is at most $64$. It is also possible to obtain more precise information
about the odd member of $\cds NG$, but we shall not need to do so.

The following result extracts from Theorem~6.1 the information
we need to prove Theorem~C.
\medskip

\iitem{(6.2) THEOREM.}~~{\sl Let $K \sbs N \nor G$, where
$K \nor G$ and $G$ is $p$-solvable for some prime $p$. Assume that
$K$ has a normal Sylow $p$-subgroup $E$ and that
$\Oh pK = K$. Suppose that $|\cds NG| = 3$ and assume that
$a = \min(\cds NG)$ lies in $\cds EG$. Then $N/E$ is abelian
and $E' \sbs \zent K$, so that $E$ is metabelian.}
\medskip

\iitem{Proof.}~~Suppose $\chi \in \irr G$ and $\chi(1) = a$, and let
$\alpha$ be an irreducible constituent of $\chi_E$. Assuming that
$\alpha$ is nonprincipal, we work to show that $N/E$ is abelian and
that $E' \sbs \zent K$. 

Observe that the stabilizer of $\alpha$ in $G$ has index dividing
$\chi(1) = a$, and thus by Lemma~5.1, this stabilizer contains $N$,
and in particular $\alpha$ is invariant in $K$. Since $E \in \syl pK$,
it follows from Corollary~6.28 of \ref\book\ that $\alpha$ extends to
some character $\beta \in \irr K$, and we see that the determinantal
order $o(\beta)$ cannot be divisible by $p$ since $\Oh pK = K$. It follows
that $o(\alpha)$ is not divisible by $p$, and thus $\alpha$ cannot be
linear since we are assuming that it is nonprincipal. It follows that
$\alpha(1) > 1$, and thus, since the $p$-group $E$ is an M-group,
$\alpha$ is induced from a proper subgroup of $E$. We conclude that
$\alpha$ reduces upon restriction to some maximal subgroup of $E$,
and thus $\alpha_{E'}$ is reducible.

Now consider a $G$-chief series containing $E'$ and $E$ among its
terms. There clearly must be some chief factor $X/Y$ with
$E' \sbs Y < X \sbs E$ such that $\alpha_Y$ reduces, but $\alpha_X$
is irreducible. Then $X/Y$ is a reducing $p$-section for $\chi$
in $N$, and we argue that $X/Y$ cannot be centralized by $N$. In
fact, we claim that $K$ acts nontrivially on $X/Y$. To see why
this is so, observe that $K/E$ acts coprimely on the abelian group
$E/E'$, and thus by Fitting's theorem,
$E/E' = [(E/E'),(K/E)] \times \cent{E/E'}{K/E}$. However, 
$[E,K] = E$ since $E$ is a normal Sylow $p$-subgroup of $K$ and
$\Oh pK = K$, and thus the first factor in the Fitting decomposition is
the whole group $E/E'$. It follows that $\cent{E/E'}{K/E}$ is trivial,
and thus $K/E$ has no nontrivial fixed points on $X/E'$. We conclude
(since this is a coprime action) that $K/E$ has no nontrivial fixed points
on $X/Y$, and in particular, the action of $K$ on $X/Y$ is nontrivial, as
claimed.

Now $X/Y$ is a reducing $p$-section for $\chi$ in $N$, but it is not
central in $N$, and thus Theorem~6.1 applies. In particular, by 6.1(d),
we see that $N/X$ is abelian, and thus $N/E$ is abelian, as required.
Also, we know that $p = 2$ and that two of the three members of
$\cds NG$ are powers of $2$, while the remaining degree is odd.

Next, we consider a $p$-complement $Q$ in $K$, and we show that it
acts trivially on $E'$. Since $(|Q|,|E'|) = 1$ and $E'$ is a $p$-group, it
suffices to show that $Q$ acts trivially on $E'/E''$, and in fact, it suffices
to consider a $G$-chief series having $E'$ and $E''$ as terms, and to
show that $Q$ centralizes each chief factor $U/V$, where
$E'' \sbs V < U \sbs E''$. We need to prove, in other words, that
$[U,Q] \sbs V$.

It suffices to show that $[U,Q] \sbs \ker\psi$ for every character
$\psi \in \irr G$ such that $V \sbs \ker\psi$. Of course, there is nothing
to prove unless $\psi \in \irrs NG$. Also, if $\psi(1)$ is odd, then the
irreducible constituents of $\psi_E$ must be linear since $E$ is a
$2$-group, and in this case, $[U,Q] \sbs U \sbs E' \sbs \ker\psi$, as
desired. We may assume, therefore, that $\psi(1)$ is a power of $2$.

Let $\gamma$ be an irreducible constituent of $\psi_K$, and note
that $\gamma_E$ is irreducible since $\gamma(1)$ is a power of $2$
and $|K:E|$ is odd. Also, $V \sbs \ker\gamma$ and $U/V$ is central
in $E/V$. Since $\gamma_E$ is irreducible, it follows that
$U \sbs \zent\gamma$. (In other words, $\gamma_U$ is a multiple of
a linear character.) Therefore $[U,Q] \sbs \ker\gamma$, and since
this holds for all choices of $\gamma$, we deduce that
$[U,Q] \sbs \ker\psi$, as desired. 

We have shown that $Q$ centralizes $E'$, and thus $\cent K{E'}$ is
a normal subgroup of $K$ having $p$-power index. Since
$\Oh pK = K$, however, we conclude that $E' \sbs \zent K$, and thus
$E'$ is abelian and $E$ is metabelian, as required.\qed
\medskip

There is an amusing consequence of the fact that the exceptional
case cannot occur in groups of odd order. Although this corollary will
not be used in what follows, it seems possible that some application
might eventually be found for it, and so we digress briefly to present it
here. To state our result, we use the notation $N^\infty$ to denote the
unique normal subgroup of a group $N$, minimal such that the
corresponding factor group is nilpotent. Note that $N^\infty$ is the final
term in the lower central series for $N$, and thus
$[N^\infty,N] = N^\infty$.
\medskip

\iitem{(6.3) COROLLARY.}~~{\sl Let $N \nor G$, where $G$ is
solvable, and let $a = \min(\cds NG)$. Assume either that $a$ is odd, or
that $|N|$ is odd, and write $M = N^\infty$. Then $a \not\in \cds MG$.}
\medskip

\iitem{Proof.}~~Let $\chi \in \irr G$ with $\chi(1) = a$. Assuming (as we
may) that $N > 1$, we see that $M < N$, and thus working by induction
on $|N|$, we can suppose that $M^\infty \sbs \ker\chi$. We want to prove
that $M \sbs \ker\chi$, and so we consider an irreducible constituent
$\alpha$ of $\chi_M$.

First, assume that $\alpha$ is linear. Observe that the stabilizer in $G$
of $\alpha$ has index dividing $\chi(1) = a$, and thus by Lemma~5.1,
$\alpha$ is invariant in $N$. It follows that $[M,N] \sbs \ker\alpha$, and
since this holds for all irreducible constituents of $\chi_M$, we deduce
that $[M,N] \sbs \ker\chi$. But $M = N^\infty$, and thus $[M,N] = M$ and
$M \sbs \ker\chi$, as required.

We complete the proof by deriving a contradiction in the case where
$\alpha$ is nonlinear. Since $M^\infty \sbs \ker\alpha$, we see that
$\alpha$ can be viewed as a character of the nilpotent group
$M/M^\infty$, and thus by the argument that we used in the proof of the
previous theorem, $\alpha_{M'}$ is reducible. If $\beta$ is an irreducible
constituent of $\alpha_{M'}$, choose any prime divisor $p$ of
$\alpha(1)/\beta(1)$ and let $U/M'$ be the $p$-complement of $M/M'$. It
is easy to see that $\alpha_U$ must be reducible, and thus there exists a
$G$-chief factor $X/Y$ with $U \sbs Y < X \sbs M$ such that $X/Y$ is a
reducing $p$-section for $\chi$ in $N$. Furthermore, $p \ne 2$ since we
are assuming that either $a$ or $|N|$ is odd, and thus by Theorem~6.1,
the section $X/Y$ is central in $N$.

Let $Q/M$ be the $p$-complement of the nilpotent group $N/M$, and
note that the coprime action of $Q/M$ on $M/U$ must have nontrivial
fixed points since the action of $Q/M$ on the section $X/Y$ is trivial.
But $M/U = [(M/U),(Q/M)] \times \cent{M/U}{Q/M}$ by Fitting's
theorem, and we know that the second factor is nontrivial, and we
deduce that $[(M/U),(Q/M)] < M/U$. It follows that there exists
a nontrivial chief factor $M/K$ of $G$, centralized by $Q$, and such that
$U \sbs K$.

Now let $P/M$ denote the Sylow $p$-subgroup of $N/M$. Then $P/M$
acts on the $p$-group $M/K$, and thus $[(M/K),(P/M)] < M/K$. Since
$M/K$ is a $G$-chief factor, it follows that $[(M/K),(P/M)] = 1$, and thus
$P$ centralizes $M/K$. We now know that $N = PQ$ acts trivially on
$M/K$, and thus $[M,N] \sbs K$. This is not the case, however, since
$[M,N] = M$, and this is the desired contradiction.\qed
\medskip

Corollary~6.3 clearly yields an alternative proof (when it applies)
of the fact that if $N \nor G$ and $N$ is solvable, then the
Fitting height of $h(N) \le |\cds NG|$. (See Corollary~4.4.) This
``top-down" argument, however, works only when $G$ is solvable and
$|N|$ is odd, while the ``bottom-up" proof in Section~4 works in all
cases.

Corollary~6.3 should be compared with the result of T.~R.~Berger
\ref\berg, which makes a stronger oddness assumption, but obtains a
very much stronger conclusion. In the case where $|G|$ is odd,
Berger proved that $a \not\in \cds{N'}G$, where $N \nor G$ and
$a = \min(\cds NG)$. Berger's conclusion does not remain valid,
however, under the weaker hypotheses of Corollary~6.3. If $G$ is
the semidirect product of $Q_8$ acting faithfully on a nonabelian
group $N$ of order $27$ and exponent $3$, for example, we see that
$\min(\cds NG) = 3$. In this situation, $N'$ is not contained in the
kernels of the irreducible characters of $G$ of degree $3$.
\bigskip

\iitem{7. Three degrees.}

In this section, we prove Theorem~C. Recall that the assumption is
that $N \nor G$, where $G$ is $p$-solvable for all prime divisors $p$ of
$|N|$, and $|\cds NG| = 3$. Our goal is to prove that $\dl N \le 3$.
\medskip

\iitem{Proof of Theorem~C.}~~Write $\cds NG = \{a,b,c\}$ with
$a < b < c$ and assume that $\dl N > 3$. Working by induction on
$|N|$ and following the usual argument, we consider $M \sbs N$ with
$M \nor G$. Then $\cds{(N/M)}{(G/M)} \sbs \cds NG$, and thus if
$M > 1$, the inductive hypothesis yields that $\dl{N/M} \le 3$ and
$N''' \sbs M$. We conclude that $N'''$ is the unique minimal normal
subgroup of $G$ contained in $N$, and thus $F = \fit N$ is a $p$-group,
where $p$ is the unique prime divisor of $|N'''|$. Also, by Corollary~3.3,
we know that $N$ is not a $p$-group, and thus $F < N$ and $N$ does
not have a normal $p$-complement. By Theorem~D, therefore, $p$ fails
to divide at least one of $a$, $b$ and $c$.

Suppose $M < N$ with $M \nor G$. Since $\cds MG \sbs \cds NG$, we
have $M''' = 1$, and thus $M''$ is abelian. Also, since $\dl N > 3$,
we see that $\dl{N/M''} > 2$, and thus $|\cds{(N/M'')}{(G/M'')}| > 2$ by
Theorem~5.3. It follows that $\cds{(N/M'')}{(G/M'')} = \cds NG$, and in
particular, $c \in \cds{(N/M'')}{(G/M'')}$. Since $\fit{M'} \sbs F$ is a
$p$-group, we can apply Lemma~4.6 to the subgroup $M' \nor G$ to
deduce that $M'$ is itself a $p$-group. This shows that $M' \sbs F$
whenever $M < N$ with $M \nor G$.

Now let $K = \Oh pN$ and write $E = K \cap F$, so that
$E = \fit K = \oh pK$. We claim that $E \in \syl pK$, and to prove this,
we consider separately the cases where $K < N$ and $K = N$.
If $K < N$, then $K' \sbs F$ by the result of the previous paragraph, and
so $K/E$ is abelian. Since $E = \oh pK$, it follows that $K/E$ is a
$p'$-group, and thus $E \in \syl pK$, as desired.

Now suppose that $K = N$ and let $N/M$ be a chief factor of $G$, so
that $N/M$ is a $q$-group for some prime $q \ne p$. Then $F \sbs M$,
and we know that $M' \sbs F$, so that $M/F$ is abelian. Since
$F = \oh pM$, it follows that $M/F$ is a $p'$-group, and hence $N/F$
is also a $p'$-group. Since we are assuming that $K = N$, we have
$E = F$, and thus $E \in \syl pK$ in this case too.

Since $\Oh pK = K$ and $E$ is a normal Sylow $p$-subgroup of $K$,
we can apply Theorem~6.2. Because $\dl N > 3$, we see that either
$N/E$ is nonabelian, or else $E$ is not metabelian, and thus 6.2
guarantees that $a \not\in \cds EG$.

Assume until further notice that $K < N$. As we have seen, this implies
that $K' \sbs F$, and thus $K' \sbs E$ and $K/E$ is abelian. Now $K$ is
not a $p$-group since otherwise, $N$ would be a $p$-group, and we
have seen that this is not the case. But $\fit K$ is a $p$-group, and thus
by Lemma~4.6, we cannot have $c \in \cds{(N/K')}{(G/K')}$, and
therefore $c \not\in \cds{(N/E)}{(G/E)}$. We deduce that
$\cds{(N/E)}{(G/E)} \sbs \{a,b\}$, and in particular, $N/E$ is metabelian
by Theorem~5.3. It follows that $E$ is nonabelian. Also, we claim that
$K' = E$ since otherwise, $K' < E$ and $p$ divides $|K/K'|$. This cannot
happen, however, since $K = \Oh pK$. 

We have seen that $a \not\in \cds EG$. If also $a \not\in \cds KG$, then
$E' = K'' = 1$ by Theorem~5.3, and $E$ is abelian, which is a
contradiction. It follows that there exists $\chi \in \irr G$ with $\chi(1) = a$,
such that $K \not\sbs \ker\chi$. Since we know that $E \sbs \ker\chi$ and
that $K/E$ is abelian, we see that an irreducible constituent $\lambda$
of $\chi_K$ is linear and nonprincipal. The stabilizer of $\lambda$ in
$G$ has index dividing $\chi(1) = a$, and thus $\lambda$ is invariant in
$N$ by Lemma~5.1. Since $E \sbs \ker\lambda$ and $K/E$ is a
$p'$-group, while $N/K$ is a $p$-group, it follows that $\lambda$
extends to $N$. An extension of $\lambda$ to $N$, however, is a linear
character whose order is a multiple of $o(\lambda)$. This order
has a nontrivial $p'$-part, and we conclude that $\Oh{p'}N < N$.

Let $L = \Oh{p'}N$ and observe that $F \sbs L$ since $F$ is a
$p$-group. Since $L < N$, we know by an earlier argument that
$L' \sbs F$, and thus $L/F$ is abelian. But $\Oh{p'}L = L$, and thus
$L/F$ must be a $p$-group, and we deduce that $L$ is a $p$-group,
and hence $L = F$ and $|N/F|$ is not divisible by $p$. It follows that
$FK = N$ and $N/E = (F/E) \times (K/E)$, and we draw a number of
conclusions. First, we observe that $K$ acts trivially on $F/E$. Also, we
know that $K/E$ is abelian but that $N/E$ is nonabelian, and so we
deduce that $F/E$ is nonabelian. Finally, since $K/E$ is abelian and
$F/E$ is a $p$-group, we see that $\cd{N/E}$ consists entirely of
$p$-powers.

Now fix $f \in \{a,b,c\}$ such that $f$ is not divisible by $p$ and suppose,
for the moment, that $f \not\in \cds EG$. (We know, for example, that
this would be the situation if $f = a$.) Let $\chi \in \irr G$ have degree $f$,
and note that $E \sbs \ker\chi$, and so the irreducible constituents of
$\chi_N$ have degrees lying in $\cd{N/E}$. The degrees of these
constituents, therefore, are $p$-powers dividing $f$, and hence
$\chi_N$ is a sum of linear characters, and $N' \sbs \ker\chi$. This
shows that $f \not\in \cds{N'}G$, and thus $|\cds{N'}G| < 3$ and $N'$ is
metabelian by Theorem~5.3. We are assuming that $\dl N > 3$,
however, and this contradiction shows that $f \in \cds EG$. In particular,
$f \ne a$, and thus $p$ divides $a$.

Again let $\chi \in \irr G$ with $\chi(1) = f$, and note that $\chi_E$ must
have linear constituents since $E$ is a $p$-group and $p$ does not
divide $f$. By the previous paragraph, we can choose $\chi$ with
$E \not\sbs \ker\chi$, and we let $\mu$ be an irreducible constituent of
$\chi_E$, so that $\mu$ is nonprincipal and linear. Let $T$ be the
stabilizer of $\mu$ in $G$ and note that $|G:T|$ divides $\chi(1) = f$,
so that $p$ does not divide $|G:T|$.

We claim now that $\mu$ extends to $T$. Since $E$ is a $p$-group, it
suffices by Theorem~6.26 of \ref\book\ to show that $\mu$ extends to
$P$, where $P \in \syl pT$. To produce an extension of $\mu$ to $P$,
we consider the Clifford correspondent $\xi \in \irr T$ of $\chi$ with
respect to $\mu$, so that $\xi_E$ is a multiple of $\mu$ and
$\xi^G = \chi$. It follows that $p$ does not divide $\xi(1)$, and thus some
irreducible constituent of $\xi_P$ has degree not divisible by $p$, and
hence is linear. This linear character is an extension of $\mu$ to $P$, as
desired, and it follows that $\mu$ extends to $T$, as claimed.

Let $\nu \in \irr T$ be an extension of $\mu$. Then $\nu^G$ is irreducible
by the Clifford correspondence, and thus $\nu^G \in \irrs NG$ since
$\nu$ extends $\mu \ne 1_E$. It follows that
$|G:T| = \nu^G(1) \in \cds NG = \{a,b,c\}$. Also, we know that $|G:T|$ is
not divisible by $p$, and thus $|G:T| \ne a$. Since $F \nor G$ is
a $p$-group, we see that $F \sbs T$, and thus $T/E$ contains the
nonabelian subgroup $F/E$. Therefore, $T/E$ is nonabelian, and
has a nonlinear irreducible character $\beta$, which we can view as a
character of $T$. (In fact, since $F/E$ is a nonabelian normal
$p$-subgroup of $T/E$, we can choose $\beta$ so that $\beta(1)$ is
divisible by $p$.) Now $\beta\nu \in \irr T$ lies over $\mu$, and thus
$(\beta\nu)^G$ is irreducible and lies in $\irrs NG$. The degree of this
character is $\beta(1)|G:T|$, which exceeds $|G:T|$, and we deduce
that $|G:T| \ne c$. The only remaining possibility is that $|G:T| = b$, and
hence $p$ does not divide $b$. Also, we must have $\beta(1)|G:T| = c$,
and thus $\beta(1) = c/b$, which is therefore an integer. Furthermore,
since $\beta$ was an arbitrary nonlinear irreducible character of $T/E$,
we deduce that $\cd{T/E} = \{1,c/b\}$, and hence $T/E$ is metabelian.
(We mention that since $\beta$ could have been chosen to have degree
divisible by $p$, we can conclude that $p$ divides $c/b$ and thus $p$
divides $c$ and $f = b$. We will not need this information, however.)

Now $F \sbs T$ and $T/E$ is metabelian, while $F/E$ is nonabelian.
It follows that $F/E \not\sbs (T/E)'$, and thus there exists a linear
character $\tau \in \irr T$ with $E \sbs \ker\tau$ and such that
$F \not\sbs \ker\tau$. It follows that $\tau^G$ is a character of degree
$|G:T| = b$, and all of its irreducible constituents lie in $\irrs NG$. There
are just two possibilities, therefore: either $\tau^G$ is irreducible, or else
it is a sum of irreducible characters of degree $a$. If the latter situation
occurs, however, then $a$ would divide $\tau^G(1) = b$, and this is
impossible since $p$ divides $a$ but $p$ does not divide $b$.
We conclude, therefore, that $\tau^G$ is irreducible.

Now write $\sigma = \tau_F$, and recall that $\sigma$ is nonprincipal
by the choice of $\tau$, and also $E \sbs \ker\sigma$. Let $S$ be the
stabilizer of $\sigma$ in $G$ and note that $S \sps T$ since $\sigma$
extends to $T$. Since $|G:T|$ is not divisible by $p$, we see that
$T$ contains a full Sylow $p$-subgroup of $S$, and thus $\sigma$
extends to this Sylow subgroup. Because $F$ is a $p$-group, it follows
by Theorem~6.26 of \ref\book\ that $\sigma$ extends to some character
$\delta \in \irr S$. Then $\delta^G$ is irreducible and has degree $|G:S|$.
Also, $\delta^G$ lies in $\irrs NG$ since $\sigma$ is nonprincipal, and it
follows that $|G:S| \in \{a,b,c\}$. Also, $|G:S|$ divides $|G:T| = b$ and
we have observed that $a$ does not divide $b$ since $p$ divides $a$.
We deduce that $|G:S| = b$, and thus $S = T$.

Recall that $K$ centralizes $F/E$, and thus $K$ stabilizes $\sigma$
since $E \sbs \ker\sigma$. Thus $K \sbs S = T$, and $K$ stabilizes
$\mu$. Since $\mu$ is a linear character of the $p$-group $E$ and
$p$ does not divide $|K:E|$, we deduce that $\mu$ extends to $K$.
Since $\mu$ is nonprincipal and has $p$-power order, its extension to
$K$ has order divisible by $p$, and this contradicts the fact that
$K = \Oh pK$. This contradiction arose from the assumption that
$K < N$, and so now we can now assume that $K = N$ and $E = F$.

Since $a \not\in \cds FG$, we see that $\cds FG \sbs \{b,c\}$. Thus $F$
has derived length at most $2$, and we conclude that $N/F$ must be
nonabelian. As before, we let $N/M$ be a chief factor of $G$ and we
write $q$ to denote the unique prime divisor of $|N/M|$, so that
$q \ne p$ and $F \sbs M$. Then $M/F$ is abelian, and all Sylow
subgroups of $N/F$ except possibly the Sylow $q$-subgroups are
abelian. If there also exists a $G$-chief factor $N/L$ that is not a
$q$-group, it follows that all Sylow subgroups of $N/F$ are abelian.
In this case, $ML = F$ and since $M/F$ and $L/F$ are abelian, we
deduce that $N/F$ is nilpotent, and thus it is abelian. We know that this
is not the case, however, and so we conclude that $L$ does not exist
and that $N = \Oh rN$ for all primes $r \ne q$.

Let $R/F$ be the normal $q$-complement of $M/F$, and suppose that
$a \not\in \cds RG$. Then $R'' = 1$, so that $R'$ is abelian, and hence
$\dl{N/R'} > 2$. We therefore have $|\cds{(N/R')}{(G/R')}| = 3$, and thus
$c \in \cds{(N/R')}{(G/R')}$. Lemma~4.6 now applies, and we conclude
that $R$ is a $p$-group, and so $R = F$ in this case.

If $R > F$, therefore, there exists $\psi \in \irr G$ with $\psi(1) = a$ and
$R \not\sbs \ker\psi$. But $F \sbs \ker\psi$ and $R/F$ is abelian, and
thus $\psi_R$ has a nontrivial linear constituent $\lambda$ with
multiplicative order not divisible by $q$. Since the stabilizer of
$\lambda$ in $G$ has index dividing $a$, we know by Lemma~5.1 that
this stabilizer must contain $N$, and it follows that $\lambda$ extends
to $N$ since $N/R$ is a $q$-group. Therefore, $\Oh rN < N$ for
all primes $r$ dividing $o(\lambda)$. This is not the case, however, and
we deduce that $R = F$, so that $N/F$ is a nonabelian $q$-group.

Now choose $e \in \{a,b,c\}$ having the smallest possible $q$-part.
Since $F'$ is abelian, we see that $\dl{G/F'} > 2$, and therefore
$\cds{(N/F')}{(G/F')} = \{a,b,c\}$. It follows that there exists
$\chi \in \irrs NG$ with $\chi(1) = e$, and such that $F' \sbs \ker\chi$.
For every such choice of $\chi$, Theorem~3.1 guarantees that
$N' \sbs \ker\chi$. (Observe that Theorem~3.1 applies because
$F/F'$ is an abelian normal $q$-complement in $N/F'$.) As
$F \sbs N' \sbs \ker\chi$, we deduce that $e \in \cds{(N/F)}{(G/F)}$,
and by Corollary~4.3, therefore, we have $e \ne c$.

We showed in the previous paragraph that if $\chi(1) = e$ and
$F' \sbs \ker\chi$, where $\chi \in \irr G$, then $N' \sbs \ker\chi$. This
cannot happen for every character $\chi \in \irr G$ of degree $e$,
however, or else $e \not\in \cds{N'}G$, and thus $N'$ is metabelian,
which is not the case. It follows that $e \in \cds{F'}G$, and in particular,
$e \ne a$ since we know that $a \not\in \cds FG$. We deduce that
$e = b$, and thus the $q$-part of $a$ exceeds the $q$-part of $b$, and
$a$ does not divide $b$. Also, $b = e \in \cds{F'}G$, and hence there
is some irreducible character of $G$ of degree $b$ that has nonlinear
constituents upon restriction to $F$. We conclude from
this that $p$ divides $b$.

Since $N' < N$, we know that $N'' \sbs F$, and thus
$a \not\in \cds{N''}G$, and $\cds{N''}G \sbs \{b,c\}$. If $p$ divides $c$,
then since $p$ also divides $b$, Theorem~D would imply that $N'$ has
a normal $p$-complement. But the unique minimal normal subgroup
of $G$ contained in $N$ is a $p$-group, and thus the $p$-complement
of $N'$ would have to be trivial. In this situation, $N'$ is a $p$-group,
and thus $N/F$ is abelian. We know that this is not the case, however,
and we conclude that $p$ does not divide $c$.

Now let $\alpha \in \irr F$ be nonlinear. (Note that $\alpha$ exists
because we know that $F$ is nonabelian.) Since irreducible characters
of $G$ of degree $a$ have $F$ in the kernel and we know that $p$ does
not divide $c$, we see that all irreducible characters of $G$ lying over
$\alpha$ must have degree $b$. Since $N/F$ is a nonabelian $q$-group,
it follows that some irreducible character of $N$ lying over $\alpha$ has
degree divisible by $q$. (This is clear if $\alpha$ does not extend to $N$.
On the other hand, if $\alpha$ extends to $\gamma \in \irr N$, let
$\beta \in \irr{N/F}$ be nonlinear, and view $\beta \in \irr N$. Then
$\beta\gamma$ is an irreducible character of $N$ lying over $\alpha$ 
and having degree divisible by $\beta(1)$, which is divisible by $q$.)
It follows that $q$ divides $b = e$, and thus $q$ divides all three
members of $\cds NG$.

Let $T$ be the stabilizer of $\alpha$ in $G$, so that $|G:T|$ divides $b$
and all irreducible characters of $T$ lying over $\alpha$ have degree
$b/|G:T|$. Given any prime $r \ne p$, we know from Lemma~5.2 that
some irreducible character of $T$ lying over $\alpha$ has degree not
divisible by $r$, and it follows that $b/|G:T|$ must be a power of $p$.

We claim now that $T \cap N = F$. To see why this is so, write
$D = T \cap N \nor T$, and suppose that $D > F$, so that there exists
a nonprincipal linear character $\tau$ of $D/N$, which we view as a
character of $D$. Let $S$ be the stabilizer of $\tau$ in $T$ and let
$Q/F \in \syl q{S/F}$. Since $\alpha$ is invariant in $Q$,
there is an extension $\delta$ of $\alpha$ to $Q$, and we can choose
$\delta$ such that $o(\delta) = o(\alpha)$ is a $p$-power. Observe that
$\delta_D$ is the unique extension of $\alpha$ to $D$ that has
$p$-power determinantal order, and thus $\delta_D$ is invariant in $T$.
Now $\tau\delta_D$ is an irreducible character of $D$ that is invariant
in $S$, and, in fact, $S$ is the full stabilizer of $\tau\delta_D$ in $T$
because $\delta_D$ and $\tau\delta_D$ uniquely determine $\tau$.

Let $\xi \in \irr S$ lie over $\tau\delta_D$. Then $\xi^T$ is irreducible
by the Clifford correspondence, and $\xi^T$  lies over $\alpha$. It follows
that $\xi(1)|T:S|  = \xi^T(1) = b/|G:T|$ is a power of $p$, and in particular,
$\xi(1)$ is not divisible by $q$. There is thus an irreducible constituent
$\eta$ of $\xi_Q$ such that $q$ does not divide $\eta(1)$. But $Q/F$ is a
$q$-group, and it follows that $\eta_F$ is irreducible, and hence
$\eta_F = \alpha$. Since also $\delta \in \irr Q$ and $\delta_F = \alpha$,
we must have $\eta = \nu\delta$, where $\nu$ is some linear character
of $Q$ with $F \sbs \ker\nu$. Also, $\eta$ lies over $\tau\delta_D$, and
so $\nu_D\delta_D = \eta_D = \tau\delta_D$, and hence $\nu_D = \tau$
and $\tau$ extends to $Q$. Applying Theorem~6.26 of \ref\book\ in the
group $S/F$, we deduce that $\tau$ extends to $S$, and we let
$\sigma \in \irr S$ be an extension of $\tau$. Thus $\sigma$ is linear,
and since $\tau$ is nonprincipal, all irreducible constituents of
$\sigma^G$ lie in $\irrs NG$. 

Recall that $\xi(1)|T:S| = b/|G:T|$. We have
$\sigma(1) = 1 < \alpha(1) \le \xi(1)$, and thus
$$
\sigma^G(1) < \xi^G(1) = \xi(1)|G:S| = \xi(1)|T:S||G:T| = b \,,
$$
and it follows that all irreducible constituents of $\sigma^G$ have
degree $a$. Thus $a$ divides $\sigma^G(1) = |G:S|$, which is a divisor
of $b$. This is a contradiction, however, since  $a$ does not divide $b$,
and we conclude that $T \cap N = F$, as claimed.

We have now shown that the stabilizer in $N$ of every nonlinear
irreducible character $\alpha$ of $F$ is exactly $F$, and thus if we
let $U \in \syl qN$, we know that the stabilizer in $U$ of $\alpha$ is trivial,
and this holds for all nonlinear irreducible characters $\alpha$ of $F$.
Recall that $F' > 1$ and let $\epsilon$ be an arbitrary nonprincipal linear
character of $F'$. If $V$ is the stabilizer in $U$ of $\epsilon$, then
$V$ permutes the irreducible constituents of $\epsilon^F$, which has
$p$-power degree. It follows that $V$ stabilizes some irreducible
constituent of $\epsilon^F$, and we note that such a constituent must
be nonlinear since $F'$ is not contained in its kernel. It follows that
$V = 1$, and thus the action of $U$ on the (nontrivial) group of linear
characters of $F$ is Frobenius. Now $U \cong N/F$, and so $U$ is a
nonabelian $q$-group, which, as we have just seen, is a Frobenius
complement. We deduce that $q = 2$ and that $N/F$ is a generalized
quaternion group.

Now let $C = \cent G{N/F}$ so that $C \nor G$ and $G/C$ is isomorphic
to a subgroup of the automorphism group of the generalized quaternion
group $N/F$. If $|N/F| > 8$, then $NC/C$, which corresponds to the inner
automorphisms of $N/F$, is nonabelian, and thus
$|\cds{(NC/C)}{(G/C)}| \ge 2$ and thus at least two of $a$, $b$ and $c$
are degrees of irreducible characters of $G/C$. But $G/C$ is a
$2$-group in this case, and so its order is not divisible by $p$. Since $b$
is divisible by $p$, we conclude that $\irrs NG$ contains an irreducible
character of degree $c$ with $C$ in its kernel. But $F \sbs C$,
and this is a contradiction by Corollary~4.3. 

We are left with the case where $N/F \cong Q_8$, and thus $G/C$
is isomorphic to a subgroup of the symmetric group $S_4$ that contains
the Klein subgroup (corresponding to $NC/C$).  If $G/C$ is a $2$-group,
it has a linear character such that $NC/C$ is not in the kernel, and
otherwise, it has an irreducible character of degree $3$ with $NC/C$
not in the kernel. It follows that $\cds NG$ contains either $1$ or $3$.
We have seen, however, that each member of $\cds NG$ is divisible
by $q = 2$, and this is our final contradiction.\qed
\vfil\eject

\centerline{REFERENCES}
\frenchspacing\parindent = 0pt
\bigskip

1.~~T. R. Berger, Characters and derived length in groups of odd 
order. J. of Algebra, {\bf 39} (1976) 199--207
\medskip

2.~~Y. G. Berkovich, Degrees of irreducible characters and normal
$p$-complements. Proc. Amer. Math. Soc. {\bf 106} (1989) 33--34.
\medskip

3.~~S. Garrison, On groups with a small number of character degrees.
Ph.D. Thesis, University of Wisconsin, Madison, 1973.
\medskip

4.~~D. Gluck, Bounding the number of character degrees of a
solvable group. J. London Math. Soc. (2) {\bf 31} (1985) 457--462.

5.~~I. M. Isaacs, Character degrees and derived length of a solvable
group. Canad. J. of Math. {\bf 27} (1975) 146--151.
\medskip

6.~~I.M. Isaacs, Character theory of finite groups, Academic Press,
New York, 1976.

\bye